%
\input amstex
\documentstyle{amsppt}
%
%
%
%
%
\def\R{\Bbb R }

\def\P{\Bbb P }

\def\Z{\Bbb Z }
\def\C{\Bbb C }

%
%
%

%
%

\define\nl{\hfil\newline}
\define\a{\alpha}
\define\la{\lambda}
\redefine\d{\delta}

\define\w{\omega}

\redefine\i{{\,\roman i\,}}
\define\mapleft#1{\smash{\mathop{\longleftarrow}\limits^{#1}}}
\define\mapright#1{\smash{\mathop{\longrightarrow}\limits^{#1}}}
\define\mapdown#1{\Big\downarrow\rlap{
   $\vcenter{\hbox{$\scriptstyle#1$}}$}}
\define\mapup#1{\Big\uparrow\rlap{
   $\vcenter{\hbox{$\scriptstyle#1$}}$}}
\define\im{\text{Im\kern1.0pt }}
\define\re{\text{Re\kern1.0pt }}

\define\pbar{\overline{\partial}}
\define\Tmp#1{T^{(m)}_{#1}}

\define\Cim{C^{\infty}(M)}

\define\ghm{\Gamma_{hol}(M,L^{m})}
\define\ghmo{\Gamma_{hol}(M,L^{m_0})}

\define\gh{\Gamma_{hol}(M,L)}

\define\gul{\Gamma_{\infty}(M,L)}

\define\Lp{{\roman L}^2(M,L)}

\define\Lqv{{\roman L}^2(Q,\nu)}

\define\rmm{\Pi^{(m)}}

\define\Qm#1{\widehat{Q}_{#1}^{(m)}}

\define\Pm#1{\widehat{P}_{#1}^{(m)}}

\define\Tfm{T_f^{(m)}}
\define\Tgm{T_g^{(m)}}
\define\Tfgm{T_{\{f,g\}}^{(m)}}
\define\Tma#1{T_{#1}^{(m)}}
\define\Hm{\Cal H^{(m)}}
\define\Hc{{\Cal H}}
\define\phm{\varPhi^{(m)}}
\define\phtm{\tilde\varPhi^{(m)}}

\define\zb{\overline{z}}
\define\wb{\overline{w}}
\redefine\d{\partial}
\define\db{\overline{\partial}}
\define\ze{\zeta}
\define\zeb{\overline{\zeta}}

\define\pfz#1{\frac {d#1}{dz}}

\define\pfzb#1{\frac {d#1}{d\overline{z}}}
\define\Pfz#1{\frac {\partial #1}{\partial z}}
\define\Pfzb#1{\frac {\partial#1}{\partial\overline{z}}}

\define\End{\text{End}}

\define\BT{Berezin-Toeplitz}
%
%
\magnification=1200
\vsize=21.5truecm
\hsize=16truecm
\hoffset=0.5cm\voffset=0.2cm
\baselineskip=15pt plus 0.2pt
\parskip=8pt
\NoBlackBoxes
\TagsOnRight
{\parskip=0pt
\hfill Mannheimer Manuskripte 203

\hfill q-alg/9601016

\hfill Rev. 28.3.96
}
\vskip 1cm
\topmatter
\title
Berezin-Toeplitz Quantization of compact
K\"ahler manifolds
\endtitle
\rightheadtext{Berezin-Toeplitz Quantization}
\leftheadtext{M. Schlichenmaier}
\author Martin Schlichenmaier
\endauthor
\address
Department of Mathematics and Computer Science,
University of Mannheim
D-68131 Mannheim, Germany
\endaddress
\email
schlichenmaier\@math.uni-mannheim.de
\endemail
\date January 96
\enddate
\abstract
In this lecture results are reviewed obtained by the author
together  with Martin Bordemann  and Eckhard Meinrenken on the
Berezin-Toeplitz
quantization of compact K\"ahler manifolds. Using global
Toeplitz operators, approximation results for the
quantum operators are shown.
{}From them it follows that the
quantum operators have the correct classical limit.
A star product deformation of the Poisson algebra is constructed.
\endabstract
\endtopmatter
{\sl  Invited lecture  at the XIV${}^{th}$  workshop
on geometric methods in physics,  Bia\l owie\D za, Poland,
July 9-15, 1995
}
%
%
%
\newcount\refCount
\def\newref#1 {\advance\refCount by 1
\expandafter\edef\csname#1\endcsname{\the\refCount}}
\widestnumber\key{99}
\newref rAAGM 
\newref rADo 
\newref rBFFLS 
\newref rBeC 
\newref rBeTQ  
\newref rBeCo 
\newref rBHSS 
\newref rBMS 
\newref rBLU 
\newref rBGTo
\newref rCGR
\newref rCaIe 
\newref rCoDe 
\newref rCoXi 
\newref rDeLe  
\newref rDeWL
\newref rFedTQ 
\newref rFlSt  
\newref rGHPA
\newref rGuCT
\newref rHoePD
\newref rKaMa  
\newref rKaMaP  
\newref rKlLeQr
\newref rMor  
\newref rMorOr  
\newref rOdCs
\newref rOMY   
\newref rRiDQ  
\newref rSCHLRS
\newref rSnGQ   
\newref rTuQM
\newref rUnUp 
\newref rWeFe  
\newref rWoGQ
%
%

%
%
%
%
\document
\head
1. Introduction\hfill\hbox{ }
\endhead

Let me start with some mathematical aspects of quantization.
 As a mathematician,
especially as an algebraic geometer,  I find the following concepts
very fascinating.
Dear reader if you are a physicist or a fellow mathematician working
in a different field (e.g\. in measure theory) you will probably
prefer other aspects of the quantization. So please excuse if these
other important concepts are not covered here.

The  arena of classical mechanics is as follows. One starts with
a phase space  $M$, which locally should represent  position and momentum.
We assume $M$ to be a differentiable manifold.
The physical observables are functions on $M$. One needs a
symplectic form $\w$, a non-degenerate antisymmetric closed
2-form, which roughly speaking opens the possibility to introduce
dynamics.
This form defines a Poisson structure on $M$ in the following
way. One assign to every function $f$ its Hamiltonian vector field
$X_f$ via
$$f\in\Cim \ \mapsto\  X_f,\quad\text{with }\  X_f\ \text{\ defined by}\quad
\w(X_f,.)=df(.)\ .$$
A Lie algebra structure on $\ \Cim\ $ is now defined  by the product
$$\{f,g\}:=\w(X_f,X_g)\ .$$
The Lie product
fulfils  the compatibility
$$
\text{for all}\quad f,g,h\in \Cim:\quad
\{f\cdot g,h\}=f\cdot\{g,h\}+\{f,h\}\cdot g\ .
$$
This says that $\ (\Cim,\,\cdot\,,\{..\,,..\})\ $ is a Poisson algebra.
The pair $\ (M,\w)\ $ is called a symplectic manifold.
A Hamiltonian system $\ (M,\w,H)\ $ is given by fixing
a function $H\in\Cim$, the
so called Hamiltonian function.

The first part of quantization
(and only this step will be discussed here)
consists in replacing the
commutative algebra of functions
by something noncommutative.
But there is the fundamental requirement, that the classical situation
(including the Poisson structure) should be recovered again as
``limit''.
There are some methods to achieve at least partially this goal.
I do not want to give a review of these methods. Let me just
mention a few.
There is the ``canonical quantization'', the
deformation quantization using star product,
geometric quantization,
Berezin quantization using coherent states
and Berezin symbols, \BT\ quantization,
and so on. I am heading here for \BT\ quantization which has
relations to the more known geometric quantization
as introduced by Kostant and Souriau.
In the following section I will recall some necessary definitions
for the case I will consider later on.
For a systematic treatment see
\cite\rSnGQ, \cite\rWoGQ.
%
\head
2. Geometric Quantization\hfill\hbox{ }
\endhead

Here I will assume $\ (M,\w)\ $ to be a K\"ahler manifold, i.e\.
$M$ is a complex manifold and $\w$ a K\"ahler form.
This says that $\w$ is a positive, non-degenerate closed 2-form of
type  $(1,1)$.
If $\dim_\C M=n$  and $z_1,z_2,\ldots,z_n$ are local holomorphic
coordinates then it can be written as
$$\w=\i\sum_{i,j=1}^ng_{ij}(z)dz_i\wedge d\zb_j\ ,\qquad
g_{ij}\in C^\infty(M,\C)\ ,$$
where the matrix $\ (g_{ij}(z))$ is
for every $z$ a positive definite  hermitian matrix.
Obviously $(M,\w)$ is a symplectic manifold.
A further data is $\ (L,h,\nabla)\ $,  with $L$ a holomorphic line bundle,
$h$ a hermitian metric on $L$
(conjugate-linear in the first argument), and $\nabla$ a connection which is
compatible  with
the metric and the complex structure.
With respect to local holomorphic coordinates
and with respect to a local holomorphic frame of the bundle
it can be given as
$\ \nabla=\d +\d\log h+\db$.
The curvature of $L$ is defined as
$$F(X,Y)=\nabla_X\nabla_Y-\nabla_Y\nabla_X-\nabla_{[X,Y]}\ .$$
\definition{Definition}
The K\"ahler manifold $(M,\w)$ is called quantizable, if there
is such a triple $(L,h,\nabla)$ with
$$F(X,Y)=-\i\w(X,Y)\ .\tag 1$$
\enddefinition
The condition (1) is called the prequantum condition.
The bundle $\ (L,h,\nabla)\ $ is called a (pre)quantum line  bundle.
Usually we will drop $h$ and $\nabla$ in the notation.
\example{Example 1}
The flat complex space $\C^n$ with
$$
\w=\i\sum_{j=1}^ndz_j\wedge d\zb_j\ .
$$
\endexample
\example{Example 2}
The Riemann sphere, the complex projective line,
 $\P(\C)=\C\cup \{\infty\}\cong S^2$. With respect to the quasi-global
coordinate $z$  the form can be given as
$$
\w=\frac {\i}{(1+z\zb)^2}dz\wedge d\zb\ .$$
The quantum line bundle  $L$ is the hyperplane  bundle.
For the Poisson bracket one obtains
$$\{f,g\}=\i(1+z\zb)^2\left(\Pfzb f\cdot\Pfz g-\Pfz f\Pfzb g\right)\ .
$$
\endexample
\example{Example 3}
The (complex-) one dimensional torus $M$.
Up to isomorphy it can be given as
$M = \C/\Gamma_\tau$ where $\ \Gamma_\tau:=\{n+m\tau\mid n,m\in\Z\}$
is a lattice with $\ \im \tau>0$.
As K\"ahler form we take
$$
\w=\frac {\i\pi}{\im \tau}dz\wedge d\zb\ ,$$
with respect to the coordinate $z$ on the
covering space $\C$.
The corresponding quantum line bundle is the theta line bundle
of degree 1, i.e\. the bundle whose global sections are
multiples of the Riemann theta function.
\endexample
\example{Example 4}
A compact Riemann surface $M$ of genus $g\ge 2$.
Such an $M$ is the quotient of the open unit disc $\Cal E$ in
$\C$ under the fractional linear transformations of
a Fuchsian subgroup of $SU(1,1)$.
If $\ R=\pmatrix a&b\\ \overline{b} &\overline{a}\endpmatrix
\ $  with $\  |a|^2-|b|^2=1\ $ (as an element of $SU(1,1)$) then
the action is
$$z\mapsto R(z):=\frac {a z + b} {\overline{b} z + \overline{a}}\ .$$
The K\"ahler form
$$\w=\frac {2\i}{(1-z\zb)^2}dz\wedge d\zb\ $$
of $\Cal E$ is invariant under the fractional linear
transformations. Hence it defines a K\"ahler form on $M$.
The quantum bundle is the canonical bundle, i.e\. the bundle
whose local sections are the holomorphic differentials.
Its global sections can be identified with the automorphic forms
of weight $2$ with respect to the Fuchsian group.
\endexample
\example{Example 5}
The complex projective space $\P^n(\C)$.  This generalizes  Example 2.
The points in  $\P^n(\C)$ are given by their homogeneous coordinates
$\ (z_0:z_1:\ldots:z_n)\ $. In the affine chart with $z_0\ne 0$ we
take $\ w_j=z_j/z_0\ $ with $j=1,\ldots, n$ as holomorphic
coordinates. The K\"ahler form is the Fubini-Study fundamental
form
$$
\w_{FS}:=
\i\frac
{(1+|w|^2)\sum_{i=1}^ndw_i\wedge
d\wb_i-\sum_{i,j=1}^n\wb_iw_jdw_i\wedge d\wb_j} {{(1+|w|^2)}^2}
\ .$$
The quantum line bundle is the hyperplane bundle $H$, i.e\. the
line bundle whose global  holomorphic  sections can
 be identified with the linear forms in the $n+1$ variables $z_i$.
\endexample
\example{Example 6}
Projective K\"ahler submanifolds.
Let $M$ be a complex submanifold of $\P^N(\C)$ and denote
by
$\ i:M\hookrightarrow \P^N(\C)$ the  embedding,  then
the pull-back of the Fubini-Study form
$i^*(\w_{FS})=\w_M$ is a K\"ahler form on $M$ and
the pull-back of the hyperplane bundle  $i^*(H)=L$ is a quantum line
bundle for the K\"ahler manifold $(M,\w_M)$.
Note that by general results $i(M)$ is an algebraic
manifold.
\endexample

There is an important observation.
If $M$ is a compact K\"ahler manifold which is quantizable then from
the prequantum condition (1)
 we get for  the Chern form of the line bundle
the relation
$$c(L)=\frac {\i}{2\pi}F=\frac {\w}{2\pi}\ .$$
This implies that $L$ is a positive line bundle. In the terminology
of algebraic geometry it is an ample line bundle.
By the Kodaira embedding theorem  $M$ can be embedded (as  algebraic
submanifold) into projective space $\P^N(\C)$ using
a basis of the global holomorphic
  sections $s_i$ of a suitable tensor power
$L^{m_0}$ of the bundle $L$
$$z\ \mapsto\  (s_0(z):s_1(z):\ldots:s_N(z))\in \P^N(\C)\ .$$
These algebraic manifolds  can be described
as  zero sets  of homogeneous polynomials.
Note that the  dimension of the space $\ghmo$
consisting of the global holomorphic sections of $L^{m_0}$, can be determined
by the Theorem of Grothendieck-Hirzebruch-Riemann-Roch, see
\cite\rGHPA, \cite\rSCHLRS.
So even if we start with an arbitrary K\"ahler manifold  the
quantization condition will force the manifold to be an algebraic manifold
and we are in  the realm of algebraic geometry.
This should be compared with the fact that there are ``considerable more''
K\"ahler manifolds than algebraic manifolds.
This tight relation between quantization and algebraic geometry can
also be found  in the  theory of coherent states as explained
by A. Odijewicz \cite\rOdCs\ and S. Berceanu \cite\rBeC.

Here a warning is in order. With the help of the embedding into
projective space we obtain by pull-back of the Fubini-Study form
another K\"ahler form on $M$ and
by pull-back of the hyperplane bundle
another quantum bundle on $M$.
As holomorphic bundles the two bundles are the same, but in general
the K\"ahler form and the metric of the bundle and hence the
connection will be different.
Essentially, these data will only coincide
 if $M$ is a K\"ahler submanifold, or
in other words if the embedding is an isometric  K\"ahler
embedding.
The situation is very much related to Calabi's diastatic function
\cite\rCaIe, \cite{\rCGR, 2nd ref.}, see also Section 4.

Now we have to deal with the functions and how to assign
operators to them.
In geometric quantization such an assignment
is given by
$$
P:(\Cim,\{.\,,.\})\to \End(\gul,[.\,,.]),\quad
f\mapsto P_f:=-\nabla_{X_f}+\i f\cdot id\ .
$$
Here $\gul$ is the space of differentiable global sections of the
bundle $L$.
Due to the prequantum condition this is a Lie homomorphism.

Unfortunately one has too many degrees of freedom. The fields depend
locally on position and momentum. Physical
 reasons imply that they should depend
only on half of them.
Such a choice of ``half of the variables'' is called a polarization.
In general there is no unique choice of polarization.
However, for K\"ahler manifolds there is
a canonical choice of coordinates: the splitting into
holomorphic and anti-holomorphic coordinates.
To obtain a polarization we consider only sections which depend
holomorphically on the coordinates.
This is called the K\"ahler (or holomorphic) polarization.

 If we denote by
$$\Pi:\gul\to\gh,$$
the projection operator from the space of differentiable
sections onto the subspace consisting of holomorphic sections  then
the quantum operators are defined as
$$
Q:\Cim\to\End(\gh),\qquad f\mapsto
Q_f=\Pi\, P_f\,\Pi \ .
$$
This map is still a linear map. But it is not a Lie homomorphism anymore.
%
\head
3. Berezin-Toeplitz Quantization\hfill\hbox{ }
\endhead

Let the situation be as in the last section. We assume
everywhere in the following
that  $M$ is compact.
We take $\ \Omega=\frac 1{n!}\w^n\ $ as volume form  on $M$.
On the space of section $\gul$  we have the scalar product
$$
\langle\varphi,\psi\rangle:=\int_M h (\varphi,\psi)\;\Omega\  ,
\qquad
||\varphi||:=\sqrt{\langle \varphi,\varphi\rangle}\ .
\tag 2
$$
Let $\Lp$ be the  L${}^2$-completion of the space of $C^\infty$-sections
of the bundle $L$   and
$\gh$ be its finite-dimensional closed subspace of holomorphic
sections.
Again let $\ \Pi:\Lp\to\gh\ $ be the projection.
\definition{Definition}
For $f\in\Cim$ the Toeplitz operator  $T_f$
is defined to be
$$ T_f:=\Pi\, (f\cdot):\quad\gh\to\gh\ .$$
\enddefinition
In words: One multiplies the holomorphic section with the
differentiable function $f$. This yields only a differentiable section.
To obtain a holomorphic section again, we have to project it back.

The linear map
$$T:\Cim\to \End\big(\gh\big),\qquad  f\to T_f\ ,$$
will be our \BT\ quantization. It is
neither a Lie algebra homomorphism nor
an associative algebra homomorphism,
because in general
$$T_f\, T_g=\Pi\,(f\cdot)\,\Pi\,(g\cdot)\,\Pi\ne
\Pi\,(fg\cdot)\,\Pi\ .$$
{}From the point of view of Berezin's approach \cite\rBeTQ, $T_f$
is the operator with contravariant symbol $f$
(see also \cite\rUnUp).
At the end of this section I will give some more references.

Due to the compactness of $M$ this defines a map
from the commutative algebra of functions to a noncommutative
finite-dimensional (matrix) algebra.
A lot of information will get lost. To recover this
information one should consider not just the bundle $L$ alone but
all its tensor powers $L^m$
and apply all the above constructions for every $m$.
In this way one obtains a family of
matrix algebras and maps
$$
\Tma {}:\Cim\to \End\big(\ghm\big),\qquad  f\to \Tma f\ .$$
This infinite family should in some sense ``approximate'' the
algebra $\Cim$.(See \cite\rBHSS\
for a definition  of such an approximation.)

For the Riemann sphere $\P(\C)$ we
 obtain  with the help of an integral kernel
the following explicit expression for
the Toeplitz operator
$$(\Tfm s)(z)=\frac {m+1}{2\pi}\int_{\C}\frac {(1+z\zeb)^m  f(\ze)s(\ze)}
{(1+\ze\zeb)^m}
\frac {\i d\ze\wedge d\zeb}{(1+\ze\zeb)^2}\ .
$$
Here the function $s$ is representing a holomorphic section
of $L^m$.
The Toeplitz operator in our situation  has always an integral kernel.
Let $k(m):=\dim\ghm$ and take
 an orthonormal  basis $s_i,\ i=1,\ldots, k(m)$  of the space
$\ghm$ then
$$(\Tfm s)(z)=
\int_M\sum_{i=1}^{k(m)}  h^{(m)}\big(s_i(w),f(w)s(w)\big)\cdot
s_i(z)\;\Omega(w)\ .
\tag 3
$$
These Toeplitz operators are still complicated but they are easier to handle
than the quantum operators. For compact $M$
we have the following relation
$$
Q_f^{(m)}=\i\cdot T_{f-\frac 1{2m}\Delta f}^{(m)}=
\i\left(\Tfm-\frac 1{2m}T_{\Delta f}^{(m)}\right)\ .
$$
This is a result of Tuynman \cite{\rTuQM, Thm.2.1} reinterpreted
in our context, see also
\cite\rBHSS.
Here the Laplacian $\Delta$ has to be calculated with respect to the
metric $\ g(X,Y)=\w(X,IY)\ $, where $I$ is the complex structure.
We see that
for $m\to\infty$ the quantum operator of geometric quantization
will asymptotically be equal to the quantum operator of
the \BT\ quantization.

For the following let us assume that $L$ is already very ample.
This says that its global sections will already do the
embedding.
If this is not the case we would have to  start  with a certain
$m_0$-tensor power of $L$ and the form $m_0\,\w$.
The following three theorems were obtained in
joint work with Martin Bordemann and Eckhard Meinrenken \cite\rBMS.
\proclaim{Theorem 1}
For every  $\ f\in \Cim\ $  there is  some $C>0$ such that
$$||f||_\infty-\frac Cm\ \le\ ||\Tfm||\ \le\
||f||_\infty\quad
\text{as}\quad m\to\infty\ .$$
Here $||f||_\infty$ is the sup-norm of $\ f\ $ on $M$ and
$||\Tfm||$ is the operator norm on $\ghm$.
In particular, we have
$ \ \lim_{m\to\infty}||\Tfm||=||f||_{\infty}$.
\endproclaim
\proclaim{Theorem 2}
For every  $f,g\in \Cim\ $ we have
$$
||m\i[\Tfm,\Tgm]-\Tfgm||\quad=\quad O(\frac 1m)\quad
\text{as}\quad m\to\infty
\ .$$
\endproclaim
The proofs can be found in the above mentioned article \cite\rBMS.
I will give some ideas of them in the next section.

These  theorems give two  approximating sequences
of maps
$$(\Cim,||..||_\infty)\to (\frak g\frak l (n,\C), ||..||_m:=\frac 1m
||..||)\qquad
f\mapsto \i m \Tma f,\quad
f\mapsto m Q_f^{(m)}\ .$$
Restricted to real valued functions the maps take values in
$\frak u(k)$, for $k=\dim\ghm$. These families of maps are only linear
maps, not Lie homomorphism with respect to the
Poisson bracket. But by Theorem 1 they are
nontrivial and by Theorem 2 they are
approximatively Lie homomorphisms.
So every Poisson algebra of a K\"ahler manifold is a
$\frak u(k),\ k\to\infty$ limit.
This was a conjecture in \cite\rBHSS\ and  our starting point
was the aim to prove
this conjecture.
In \cite\rBMS\ also a Egorov type theorem is presented.

If one puts $\hbar=\frac {1}m$ in Theorem 2 one can rewrite
it as
$$\lim_{\hbar\to 0}
||\frac{\i}{\hbar}[T^{(1/\hbar)}_f,T^{(1/\hbar)}_g]-
T^{(1/\hbar)}_{\{f,g\}}||
              =  0\ .
$$
One should compare this with the definition of a
star product deformation of $\Cim$ (see \cite\rBFFLS, \cite\rWeFe)
based on the deformation theory of algebras as developed by Gerstenhaber.
Because there are different variants let me recall the
definition we are using.

Let $\Cal A=\Cim[[\hbar]]$ be the algebra of formal power series in the
variable $\hbar$ over the algebra $\Cim$. A product $*$ on $\Cal A$ is
called a (formal) star product if it is an
associative $\C[[\hbar]]$-linear product such that
\roster
\item
$\Cal A/\hbar\Cal A\cong\Cim$, i.e\.\quad $f*g \bmod \hbar=f\cdot g$,
\item
$\dfrac 1\hbar(f*g-g*f)\bmod \hbar = -\i\{f,g\}$.
\endroster
Note that
$\ f*g=\sum\limits_{i=0}^\infty C_i(f,g)\hbar^i\ $ with $\C$-bilinear maps
$C_i:\Cim\times\Cim\to\Cim$.
With this we calculate
$$C_0(f,g)=f\cdot g,\quad\text{and}\quad
C_1(f,g)-C_1(g,f)=-\i\{f,g\}\ .\tag 4$$
\proclaim{Theorem 3}
There exists a unique (formal) star product on $\Cim$
$$f * g:=\sum_{j=0}^\infty \hbar^j C_j(f,g),\quad C_j(f,g)\in
C^\infty(M),
\tag 5$$
in such a way that for  $f,g\in\Cim$ and for every $N$  we have
$$||T_{f}^{(m)}T_{g}^{(m)}-\sum_{0\le j<N}\left(\frac 1m\right)^j
T_{C_j(f,g)}^{(m)}||=K_N(f,g) \left(\frac 1m\right)^N\ \tag 6$$
for  $m\to\infty$, with suitable constants $K_N(f,g)$.
\endproclaim
We do not say anything about the convergence of the series (5).
Hence we do not claim to obtain a ``strict deformation quantization''
as introduced by Rieffel \cite\rRiDQ.

We obtain a star product deformation not just by
cohomological techniques as \cite\rDeWL
\ but a geometrically induced one.
There are other geometric constructions of a star product deformation
for  Poisson algebras. An important one is  given by
 Fedosov \cite\rFedTQ. {}\footnote
{See also the Bourbaki expos\'e  by Weinstein \cite\rWeFe\ 
and the review by Flato and Sternheimer \cite{\rFlSt}.}
(See  Omori, Maeda, Yoshioka \cite\rOMY\ 
and Karasev, Maslow \cite{\rKaMaP} for  related ones.)
As pointed out by Deligne \cite\rDeLe\ it would be interesting
to examine the relations between the two different approaches.

Here it is not the place and in fact I am not the expert to give
a complete list of references on the \BT\ quantization. So let me just
quote few of them.
\BT\ quantization was mainly examined for certain complex
symmetric domains. For older work see besides Berezin \cite\rBeTQ\ also
Berger-Coburn \cite\rBeCo. Similar results as stated in
Theorem 1 and Theorem 2 were recently obtained in these
cases. To give a few names:
Klimek-Lesniewski \cite\rKlLeQr, Borthwick-Lesniewski-Upmeier
\cite\rBLU, Coburn \cite\rCoDe, ...
As I will explain in Section  4 the techniques in these cases are
very different from ours.
They will not work in the case of a general  K\"ahler
manifold. On the contrary, our methods are closely related to
the compactness. So the results are at two different
edges of the theory.
Let me add that the case of compact Riemann surfaces of arbitrary
genus has been proven by the ``classical techniques''
(\cite{\rKlLeQr, 2nd ref.} for $g\ge 2$ and \cite\rBHSS\ for $g=1$).
In some cases the relation to star product deformations have been
studied \cite\rCoXi.

Closely related to the \BT\ quantization is the
quantization via
Berezin's coherent states  using the Berezin symbols \cite\rBeTQ
\ in the formulation of Cahen, Gutt and Rawnsley
\cite{\rCGR}.
This technique was also used to define star products.
See also the construction of star products by
Moreno and Ortega-Navarro \cite{\rMorOr}, \cite{\rMor}.
For the idea of relating asymptotics to a deformation
of the Poisson  bracket see Karasev and Maslov
\cite{\rKaMa}

Let me close this section with the remark that \BT\ quantization
fits into the concept of ``prime quantization'' introduced by
Ali and Doebner
\cite\rADo, \cite\rAAGM.
%
\head
4. Some remarks on the proofs\hfill\hbox{ }
\endhead

One way to prove the theorems is to represent the
sections of $L$ in a certain way, write down the
projection operator as integral operator and calculate
norms of the Toeplitz operators.
This was done by Bordemann, Hoppe, Schaller and
Schlichenmaier in \cite\rBHSS\ for the case of the
$n$-dimensional complex torus using theta functions, and
for the Riemann sphere (unpublished).
For Riemann surfaces of genus $g\ge 2$ it was done by
Klimek and Lesniewski \cite\rKlLeQr\ using
automorphic forms.
Similar techniques work for symmetric domains.
In all these cases it was important that one could
represent the sections as ordinary functions on some
simple covering of the manifold under consideration.

For general K\"ahler manifolds this does not work.
We need a different approach.
The principal idea is to group all $\Tfm $ together
to a global object.
Take $\ (U,k):=(L^*,h^{-1})\ $ the dual of the quantum line bundle,
$Q$ the unit circle bundle inside $U$ (with respect to the metric $k$)  and
$\tau: Q\to M$ the projection.
Note that for the projective space the bundle $U$ is just the tautological
bundle whose fibre over the point $z\in\P^N(\C)$ consists of
the line in $\C^{N+1}$ which is represented by $z$. In particular
the total space of $U$ without the zero section can be identified
with $\C^{N+1}\setminus\{0\}$.

Starting from the function
$\hat k(\la):=k(\la,\la)$ on $U$ we define
 $\tilde a:=\dfrac {1}{2\i}(\d-\pbar)\log \hat k$ on
$U\setminus 0$ (with respect to the complex
structure on $U$)
and restrict it to $Q$. Denote this restriction by $\alpha$.
Now $d\a=\tau^*_Q\w$ (with $d=d_Q$) and $\nu=\frac 1{2\pi}\tau^*\Omega
\wedge \a$ is a volume form
on $Q$. With respect to this form we take the L${}^2$-completion $\Lqv$
of the space of functions on $Q$.
The generalized Hardy space $\Hc$ is the closure of the
functions of $\Lqv$ which can be extended to
holomorphic functions on the whole
disc bundle.
The generalized Szeg\"o projector is the projection
$\Pi:\Lqv\to \Hc$.

 By the natural circle action
$Q$ is a $S^1$-bundle and the tensor powers of $U$ can be
viewed as associated bundles. The space $\Hc$ is preserved
by this action.
It is the (completed) direct sum
$\Hc=\sum_{m=0}^\infty \Hm$ where
 $c\in S^1$ acts on $\Hm$ as multiplication
by $c^m$.
Sections of $L^m=U^{-m}$ can be identified with functions $\phi$ on $Q$ which
satisfy the equivariance condition
$\phi(c\la)=c^m\phi(\la)$.
This identification is an isometry.
 Hence, restricted to the holomorphic objects
$$\ghm\cong\Hm\ .$$

There is the notion of Toeplitz structure
$(\Sigma,\Pi)$ as developed by
Guillemin and Boutet de Monvel
in \cite{\rBGTo}, \cite\rGuCT.
Here is not the   place to go into the details of the general
definitions.
Let me just explain what is needed here.
Here  $\Sigma$ is the symplectic submanifold
of the tangent bundle of $Q$ with the zero section
removed,
$$\Sigma=\{\;t\alpha(\lambda)\;|\;\lambda\in Q,\,t>0\ \}\ \subset\  T^*Q
\setminus 0\ ,$$
and $\Pi $ is the above projection.
A (generalized) Toeplitz operator of order $k$ is  an operator
$A:\Hc\to\Hc$ of the form
$\ A=\Pi\cdot R\cdot \Pi\ $ where $R$ is a
pseudodifferential operator
($\Psi$DO) of order $k$ on
$Q$.
The Toeplitz operators build a ring.
The (principal) symbol of $A$ is the restriction of the
principal symbol of $R$ (which lives on $T^*Q$) to $\Sigma$.
Note that $R$ is not fixed by $A$ but
Guillemin and Boutet de Monvel showed that the (principal) symbols
are well-defined and that they obey the same rules as the
symbols of   $\Psi$DOs
$$
\sigma(A_1A_2)=\sigma(A_1)\sigma(A_2),\qquad
\sigma([A_1,A_2])=\i\{\sigma(A_1),\sigma(A_2)\}_\Sigma.
\tag 7
$$
Here we use the 2-form $\omega_0=\sum_i dq_i\wedge dp_i$ on
$T^*Q$ to define the Poisson bracket there.
We are only dealing with two Toeplitz operators:
\nl
(1) The generator of the circle action
gives the  operator $D_\varphi=\dfrac 1{\i}\dfrac {\partial}
{\partial\varphi}$. It is an operator of order 1 with symbol $t$.
It operates on $\Hm$ as multiplication by $m$.
\nl
(2) For $f\in\Cim$ let $M_f$ be the multiplication operator on
$\Lqv$, i.e\.  $M_f(g)(\la):=f(\tau(\la))g(\la)$.
We set $\ T_f=\Pi\cdot M_f\cdot\Pi:\Hc\to\Hc\ $.
Because $M_f$ is constant along the fibres, $T_f$ 
commutes with the circle action.
Hence
$\ T_f=\bigoplus\limits_{m=0}^\infty\Tfm\ $,
where $\Tfm$ is the restriction of $T_f$ to $\Hm$.
After the identification of $\Hm$ with $\ghm$ we see that these $\Tfm$
are exactly the Toeplitz operators  $\Tfm$ introduced in Section 3.
In this sense we call $T_f$ also the global Toeplitz operator and
the $\Tfm$ the local Toeplitz operators.
$T_f$ is an operator of order $0$ and its symbol is just
$f$ pull-backed to $Q$ and further to $T^*Q$ (and restricted to $\Sigma$).
Let us denote by
$\ \tau^*_\Sigma:\Sigma\supseteq\tau^*Q\to Q\to M$ the composition
then we obtain for its symbol  $\sigma(T_f)=\tau^*_\Sigma(f)$.
\nl
This is the set-up more details can be found in \cite\rBMS.

\demo{Proof of Theorem 2}
Now we are able to proof Theorem 2.
The commutator
$[T_f,T_g]$ is a  Toeplitz operator of order $-1$.
Using $\ {\omega_0}_{|t\alpha(\lambda)}=-t\tau_\Sigma^*\omega\ $
for $t$ a fixed positive number, we obtain
\footnote{
Unfortunately, in \cite{\rBMS} the minus sign was missing.
This causes in Thm. 4.2 of that article also the wrong sign.}
with (7) that its
principal symbol is
$$\sigma([T_f,T_g])(t\alpha(\lambda))=\i\{\tau_\Sigma^* f,\tau_\Sigma^*g
\}_\Sigma(t\alpha(\lambda))=
-\i t^{-1}\{f,g\}_M(\tau(\lambda))\ .$$
Now consider
$$ A:=D_\varphi^2\,[T_f,T_g]+\i D_\varphi\, T_{\{f,g\}}\ .
$$
Formally this is an operator of order 1.
Using $\ \sigma(T_{\{f,g\}})=\tau^*_\Sigma \{f,g\}$ 
and $\sigma(D_\varphi)=t$ we see that its principal
symbol vanishes. Hence
it is an operator of order 0.
Now $M$ and hence $Q$ are compact manifolds.
 This implies that $A$ is a bounded
operator ($\Psi$DOs of order 0 are bounded).
It is obviously $S^1$-invariant and we can write
$A=\sum_{m=0}^\infty A^{(m)}$
where $A^{(m)}$ is the restriction of $A$ on the space $\Hm$.
For the norms we get $\ ||A^{(m)}||\le ||A||$.
But
$$
A^{(m)}=A_{|\Hm}=m^2[\Tfm,\Tgm]+\i m\Tfgm.
$$
Taking the norm bound and dividing it by $m$ we get the claim of Theorem 2.
\qed
\enddemo

\demo{Proof of Theorem 3}
This proof is a modification of the above approach.
One constructs inductively
$C_j(f,g)\in\Cim$ such that
$$A_N=D_\varphi^N T_fT_g- \sum_{j=0}^{N-1}
D_\varphi^{N-j}T_{C_j(f,g)}$$
is a zero order Toeplitz operator. Because $A_N$ is $S^1$-invariant
and it is of zero order
its principal symbol descends to a function on $M$.
Take this function to be $C_N(f,g)$. Then
$A_N-T_{C_N(f,g)}$ is of order $-1$ and
$A_{N+1}=D_\varphi(A_N-T_{C_N(f,g)})$ is of order zero.
The induction starts with
$A_0=T_fT_g$ which implies $\sigma(A_0)=\sigma(T_f)\sigma(T_g)=f\cdot g=
C_0(f,g)$.
As a zero order operator $A_N$ is bounded, hence this is true
for the component operators
$A_N^{(m)}$. We obtain
$$
||m^N\Tfm\Tgm-\sum_{j=0}^{N-1}m^{N-j}T^{(m)}_{C_j(f,g)}||\le ||A_N||\ .$$
dividing this by $m^N$ we obtain the asymptotics (6) of the theorem.
Writing this explicitly for $N=2$ we obtain for the pair $(f,g)$
$$
||m^2 \Tfm\Tgm-m^2 T^{(m)}_{f\cdot g}-m T^{(m)}_{C_1(f,g)}||
\le K\ ,
$$
 and a similar expression for the pair $(g,f)$.
By subtracting the corresponding operators,
using the triangle inequality,
dividing by $m$ and multiplying with $\i$ we obtain
$$||m\i (\Tfm\Tgm-\Tgm\Tfm)-T^{(m)}_{\i\big(C_1(f,g)-C_1(g,f)\big)}||
=O(\frac 1m)\ .$$
With Theorem 2 this yields
$\ || T^{(m)}_{\{f,g\}-\i\big(C_1(f,g)-C_1(g,f)\big)}||=O(\frac 1m)$.
But Theorem 1 says that the left hand side has as limit
$\ |\{f,g\}-\i(C_1\big(f,g)-C_1(g,f)\big)|_\infty\ $, hence
$\{f,g\}=\i(C_1(f,g)-C_1(g,f))$.
This shows equation (4). Uniqueness of the $C_N(f,g)$ follows
inductively in the same way  from (6), again using Theorem 1.
The associativity follows from the definition by
operator products.
\qed
\enddemo

Unfortunately, Theorem 1 has a rather complicated proof using
Fourier integral operators, oscillatory integrals
 and Berezin's coherent states.
(At least we have not been able to find a simpler one).
 For the special situation
of projective K\"ahler submanifolds we have a much less involved proof,
using Calabi's diastatic function.

Recall from Section 2 that a projective K\"ahler submanifold is
a K\"ahler manifold $M$ which can be embedded into  projective space
$\P^N(\C)$ (with $N$ suitable chosen) such that the K\"ahler form of $M$
coincides with the pull-back of the Fubini-Study form.
The pull-back of the tautological bundle is the dual of
the quantum bundle. We denote this bundle by $U$. On the
tautological bundle we have the standard
hermitian metric $\ k(z,w):=\langle z,w\rangle
=\bar zw\ $ in $\C^{N+1}$. By pull-back
this defines a metric on $U$. Note that in this case
the pull-back is essentially just the restriction of all objects
to the submanifold.
The Calabi (diastatic) function
\cite{\rCaIe},\cite{\rCGR, 2nd ref.}
is defined as
$$D:M\times M\to
\R_{\ge 0}\cup\{\infty\},\qquad
D(\tau(\lambda),\tau(\mu))=-\log
|k(\lambda,\mu)|^2$$
(where we have to choose $\lambda$ and $\mu$ with
$k(\lambda,\lambda)=k(\mu,\mu)=1$
representing  the points of $M$). It is well-defined,
vanishes only along the diagonal
and is strictly positive outside the diagonal.
\demo{Proof of Theorem 1 for this case}
First the easy part (which of course works in all cases).
Note that
$\ ||\Tfm||=||\rmm\,M_f^{(m)}\,\rmm||\le ||M_f^{(m)}||\ $ and
for $\varphi\ne 0$
$$\frac {{||M_f^{(m)}\varphi||}^2}{||\varphi||^2}=
\frac {\int_M h^{(m)}(f\varphi,f\varphi)\Omega}
 {\int_M h^{(m)}(\varphi,\varphi)\Omega}
=
\frac {\int_M \overline{f(z)}f(z)h^{(m)}( \varphi,\varphi)\Omega}
{\int_M h^{(m)}(\varphi,\varphi)\Omega}
\le
||f||{}_\infty^2\ .$$
Hence,
$$
||\Tfm||\le ||M_f^{(m)}||=\sup_{\varphi\ne 0}
\frac {||M_f^{(m)}\varphi||}{||\varphi||}\le ||f||_\infty .$$
To proof the first inequality, let $x_0\in M$ be a point where $|f|$ assumes
its maximum, and fix a $\lambda_0\in \tau^{-1}(x_0)$ with
$k(\lambda_0,\lambda_0)=1$.
We define a sequence
of holomorphic functions
$\phtm$ by setting
$\ \phtm(\lambda):=k(\lambda_0,\lambda)^m
\ $.
Because
$ \phtm(c\la)=c^mk(\la_0,\la)^m=c^m\phtm(\la)\ $
this defines an element $\phm$ of $\ghm$.
Note that
$$h^m(\phm,\phm)(x)
=\overline{\phtm(\la)}\phtm(\la)=
\overline{k(\la_0,\la)}^m k(\la_0,\la)^m=
\exp(-mD(x_0,x))\ $$
With
 Cauchy-Schwartz's inequality
we obtain
$$\gather
||\Tfm||\ge \frac {||\Tfm\phm||}{||\phm||}
\ge \frac {|<\phm,\Tfm\phm>|}{<\phm,\phm>}
\\=
\frac {|\int_Mf(x)h^m(\phm,\phm)(x)\Omega(x)|}
{\int_Mh^m(\phm, \phm)(x)\Omega(x)}=
\frac {|\int_Mf(x)e^{-mD(x_0,x)}\Omega(x)|}
{\int_Me^{-mD(x_0,x)}\Omega(x)}\ .
\endgather
$$
We want to consider the $\ m \to\infty\ $ limit.
The part of the integral outside a small neighbourhood 
of $x_0$ will vanish exponentially. For the rest
the stationary
phase theorem \cite\rHoePD\ allows one to
compute the asymptotics. The point $x=x_0$ is a zero of $D$ and it
is a non-degenerate critical point. Hence we obtain for the right hand
side the asymptotic
$$
\frac {|f(x_0)|+O(m^{-1})}{1+O(m^{-1})}=
|f(x_0)|+O(m^{-1}), $$
and hence
$$||\Tfm||\ge |f(x_0)|+O(m^{-1})=
||f||_\infty+O(m^{-1}) \ .\qed$$
\enddemo
%
%
%
%
%
%
%
%
{\bf Acknowledment.}
I like to thank the organizers of the conference and the audience
for the very lively atmosphere, the stimulating discussions
and the warm hospitality experienced
at the conference. My very special thanks go to
Anatol Odzijewicz   and Aleksander Strasburger.

\def\LMP{Lett\. Math\. Phys\. }
\def\CMP{Commun\. Math\. Phys\. }

\def\Izv{Math\. USSR Izv\. }

\def\DG{J\. Diff\. Geo\.}
\def\PRA{Phys\. Rev\. A}
\def\TAMS{Trans\. Amer\. Math\. Soc\.}
\def\PSPM{Proc\. Symp\. Pure Math\.}
\def\JFA{J\. Funct\. Anal\.}
\def\Adva{Advances in Math\.}
{}
\enddocument
\bye